\documentclass{article}

\usepackage{amssymb, theorem, amsmath, rotate}
\usepackage[all,2cell,dvips]{xy} \UseAllTwocells \SilentMatrices

\title{A K3 surface associated to certain
  integral matrices with integral eigenvalues}
\author{Ronald van Luijk}

\newcounter{nootje}
\newcommand\donecheck[1]{}

\newcommand\checkgone[1]{}

\newcommand\eindproof{\unskip\nobreak\hfill\hbox{\quad $\square$}\par \medskip}

\makeatletter
\def\eqalign#1{\null\,\vcenter{\openup\jot\m@th
  \ialign{\strut\hfil$\displaystyle{##}$&$\displaystyle{{}##}$\hfil
      \crcr#1\crcr}}\,}
\makeatother

\newcounter{sect}
\newcommand\mysection[1]{{\addtocounter{sect}{1}} \setcounter{theorem}{0}
                         \smallskip \noindent {\large \bf \thesect.\,\, #1} 
                         \newline \smallskip}
\newcommand\mysectiononumber[1]
                      {{\addtocounter{sect}{1}} \setcounter{theorem}{0}
                         \smallskip \noindent {\large \bf #1} 
                         \newline \smallskip}

\renewcommand\abstract[1]
    {\begin{center}
     \small
     \begin{minipage}{4in}#1\end{minipage}
     \end{center}
    }
\newcommand\myref[1]{[#1]}
\newtheorem{theorem}{Theorem}
\newtheorem{proposition}[theorem]{Proposition}
\newtheorem{definition}[theorem]{Definition}
\newtheorem{lemma}[theorem]{Lemma}
\newtheorem{corollary}[theorem]{Corollary}
\newenvironment{proof}{\noindent {\bf Proof.}}{\eindproof}
\newenvironment{proofof}{\noindent {\bf Proof of }}{\eindproof}

{\theorembodyfont{\rmfamily}

\newtheorem{remark}[theorem]{Remark}}

\renewcommand\O{\mathcal{O}}
\renewcommand\L{\mathcal{L}}
\newcommand\I{\mathcal{I}}

\newcommand\Z{\mathbb{Z}}
\newcommand\Q{\mathbb{Q}}

\newcommand\C{\mathbb{C}}

\renewcommand\P{\mathbb{P}}

\newcommand\Spec{\mathop{\rm Spec} \nolimits}

\newcommand\disc{\mathop{\rm disc} \nolimits}

\newcommand\Pic{\mathop{\rm Pic} \nolimits}
\renewcommand\deg{\mathop{\rm deg} \nolimits}
\newcommand\NS{\mathop{\rm NS} \nolimits}
\newcommand\Hom{\mathop{\rm Hom} \nolimits}
\newcommand\SL{\mathop{\rm SL} \nolimits}

\newcommand\tors{{\mathop{\rm tors} \nolimits}}

\newcommand\rk{\mathop{\rm rk} \nolimits}

\renewcommand\mod{\mathop{\rm mod} \nolimits\,}
\newcommand\Aut{\mathop{\rm Aut} \nolimits\,}

\newcommand\reg{\text{\rm reg}}
\newcommand\semi{\rtimes}

\newcommand\isom{\cong}
\newcommand\ra{\rightarrow}

\newcommand\bX{Y}
\newcommand\bXbar{\overline{Y}}
\newcommand\Xbar{\overline{X}}
\newcommand\sF{\mathcal{F}}

\newcommand\myrefart[6]{\item[{[#1]}] #2, #3, {\em #4}, #5, pp. #6.}
\newcommand\myrefbook[5]{\item[{[#1]}] #2, {\em #3}, #4, #5.}
\newcommand\nik{Ni}
\newcommand\BT{BT}
\newcommand\blv{BLV}
\newcommand\BPV{BPV}
\newcommand\bremner{Br}
\newcommand\naw{NAW}
\newcommand\shioda{Sh}
\newcommand\hag{Ha}
\newcommand\silv{Si1}
\newcommand\silvtwo{Si2}
\newcommand\tate{Ta}
\newcommand\chin{Ch}
\newcommand\artin{Ar}
\newcommand\egatwo{EGA IV(2)}
\newcommand\luijkheron{Lu}
\newcommand\sgasix{SGA 6}
\newcommand\inose{In}
\newcommand\SI{SI}

\newcommand\titleone{Introduction}
\newcommand\titletwo{Lattices and elliptic surfaces}
\newcommand\titlethree{Proof of the main theorem}
\newcommand\titlefour{The Mordell-Weil group and the N\'eron-Severi group}
\newcommand\titlefive{The surface $\bXbar$ is not Kummer}
\newcommand\titlesix{All curves on $X$ of low degree}
\newcommand\titleref{References}

\newcommand\parag{\S}

\begin{document}

\begin{center}
\Large A K3 surface associated to certain
  integral matrices with integral eigenvalues \\
\bigskip
\large Ronald van Luijk \\
Department of Mathematics 3840 \\
970 Evans Hall \\
University of California \\
Berkeley, CA 94720-3840 \\
\verb|rmluijk@math.berkeley.edu|
\end{center}

\newpage

\noindent{\bf Abstract}:
\abstract{In this article we will show that there are infinitely many
  symmetric, integral $3 \times 3$ matrices, with zeros on the
  diagonal, whose eigenvalues are all integral. We will do this by
  proving that the rational points on a certain non-Kummer, singular 
  K3 surface
  are dense. We will also compute the entire N\'eron-Severi group of
  this surface and find all low degree curves on it.
}
\bigskip

\noindent {\bf Keywords}: 
\abstract{symmetric matrices, eigenvalues, elliptic surfaces, 
K3 surfaces, N\'eron-Severi group, rational curves, Diophantine equations, 
arithmetic geometry, algebraic geometry, number theory}

\newpage

\begin{center}
\begin{tabular}{rl}
1.& \titleone \hfill \pageref{matsecone}\\
2.& \titletwo \hfill \pageref{matsectwo}\\
3.& \titlethree \hfill \pageref{matsecthree}\\
4.& \titlefour \hspace{2cm}\hfill \pageref{matsecfour}\\
5.& \titlefive \hfill \pageref{matsecfive}\\
6.& \titlesix \hfill \pageref{matsecsix}\\
  & \titleref \hfill \pageref{matsecref}\\
\end{tabular}
\end{center}
\bigskip

\frenchspacing

\mysection{\titleone}\label{matsecone}

\noindent
In the problem section of Nieuw Archief voor Wiskunde \myref{\naw}, 
F.~Beukers posed the question whether symmetric, integral $3\times 3$ matrices
\begin{equation}\label{matmatrix}
M_{a,b,c} = 
\left(
\begin{array}{ccc}
0&a&b \cr
a&0&c\cr
b&c&0\cr
\end{array}
\right)
\end{equation}
exist with integral eigenvalues and satisfying $q(a,b,c) \neq 0$, where
$q(a,b,c)$ is the polynomial
$q(a,b,c)=abc(a^2-b^2)(b^2-c^2)(c^2-a^2)$. As it is easy to find such
matrices satisfying $q(a,b,c)=0$, we will call those trivial.
R.~Vidunas and the author of this article
independently proved that the answer to this question is positive, 
see \myref{\blv}. There are in fact infinitely many nontrivial examples
of such matrices. This follows immediately from the fact that for every 
integer $t$, if we set
\begin{equation}\label{matlowdegparam}
\eqalign{
 a &= -(4t-7)(t+2)(t^2-6t+4),\cr
 b &= (5t-6)(5t^2-10t-4),\cr
 c &= (3t^2-4t+4)(t^2-4t+6),\cr
 x &= 2(3t^2-4t+4)(4t-7),\cr
 y &= (t^2-6t+4)(5t^2-10t-4),\cr
 z &= -(t+2)(5t-6)(t^2-4t+6),\cr
}
\end{equation}
then the matrix $M_{a,b,c}$ has eigenvalues $x,y$, and $z$.
This matrix is trivial if and only if we have $t\in \{-2,-1,0,1,2,4,10\}$.
For $t=3$ we get $a=125$, $b=99$, and $c=57$ with eigenvalues $190$,
$-55$, and $-135$. By a computer search, we find that this is the second 
smallest example, when ordered by
$\max(|a|,|b|,|c|)$. The smallest has $a=26$, $b=51$, and $c=114$.
In this article we will show how to find such
parametrizations. We will see that there are infinitely many and that
the one in (\ref{matlowdegparam}) has the lowest possible degree.

If the eigenvalues of the matrix $M_{a,b,c}$
are denoted by $x,y$, and $z$, then its characteristic polynomial 
can be factorized as
$$
\lambda^3 - (a^2+b^2+c^2)\lambda -2abc = (\lambda-x)(\lambda-y)(\lambda-z).
$$
Comparing coefficients, we get three homogeneous equations in
$x,y,z,a,b$, and $c$. 
Hence, geometrically we are looking for rational points on the 
$2$-dimensional complete intersection $X\subset \P^5_{\Q}$, given by
\begin{equation}\label{mateqs}
\eqalign{
x+y+z &=0, \cr
xy+yz+zx &= -a^2-b^2-c^2, \cr
xyz &= 2abc. \cr
}
\end{equation}
The points on the curves on $X$ defined by $q(a,b,c)=0$ correspond to
the trivial matrices.
Parametrizations as in (\ref{matlowdegparam}) correspond to 
curves on $X$ that are isomorphic over $\Q$ to $\P^1$. 
We will see that $X$ contains infinitely many of them, 
thereby proving the main theorem of this paper, which states the following.

\begin{theorem}\label{matmain}
The rational points on $X$ are Zariski dense.
\end{theorem}

In the next section we will recall the definition and some properties of 
lattices and elliptic surfaces in the sense of Shioda \myref{\shioda}. 
In section 3 we will prove Theorem \ref{matmain} using an elliptic fibration of
a blow-up $\bX$ of $X$. We will see that $\bX$ is a so called elliptic
K3 surface. 
The interaction between the geometry and the arithmetic of K3 surfaces is
of much interest. F. Bogomolov and Y. Tschinkel have proved that on
every elliptic K3 surface $Z$ over a number field $K$ the rational
points are potentially dense, i.e., there is a finite field
extension $L/K$, such that the $L$-points of $Z$ are dense in $Z$, see
\myref{\BT}, Thm. 1.1. 
Key in their analysis of potential density of rational points is the so 
called Picard number of a surface, an important geometric invariant.
F. Bogomolov and Y. Tschinkel have shown that if the Picard number of
a K3 surface is large enough, then the rational points are potentially
dense. On the other hand, it is not yet known if there
exist K3 surfaces with Picard number $1$
on which the rational points are not potentially dense.

After proving the main theorem, we will investigate more deeply
the geometry of $\bX$ and show in Section 4 that its Picard number
equals $20$, 
which is maximal among K3 surfaces in characteristic $0$. It is a fact
that a K3 surface with maximal Picard number is either a Kummer
surface or a double cover of a Kummer surface. These Kummer surfaces
are K3 surfaces with a special geometric structure, described in
section 5. As a consequence, their arithmetic can be described more
easily. It is therefore natural to ask if $\bX$ is a Kummer surface,
in which case $\bX$ would have had a richer structure that we could
have utilized. In Section 5 we will show that this is not the case. 

In Section 6 we will describe more of the geometry of $X$ by showing 
that $X$ contains exactly $63$ curves of degree smaller than $4$. All
points on these curves correspond to matrices that are either
trivial or not defined over $\Q$. As the degree of a parametrization 
as in (\ref{matlowdegparam}) corresponds to the degree of the curve that
it parametrizes, this shows that the one in 
(\ref{matlowdegparam}) has the lowest possible degree among
parametrizations of nontrivial matrices.

The author would like to thank Hendrik Lenstra for very helpful 
discussions. \bigskip

\mysection{\titletwo}\label{matsectwo}

We will start with the definition of a lattice. Note that for any abelian
groups $A$ and $G$, a symmetric bilinear map $A \times A
\rightarrow G$ is called {\em nondegenerate} if the induced homomorphism 
$A \rightarrow \Hom(A,G)$ is injective. Note that we do not require a
lattice to be definite, only nondegenerate.

\begin{definition}
A {\em lattice} is a free $\Z$-module $L$ of finite rank, endowed with a
symmetric, bilinear, nondegenerate map $\langle
\underline{\,\,\,\,}\,,\underline{\,\,\,\,} \rangle \colon L \times L
\rightarrow \Q$, called the {\em pairing} of the lattice. 
An {\em integral lattice} is a lattice whose pairing is $\Z$-valued. 
A lattice $L$ is called {\em even} if $\langle x,x
\rangle \in 2\Z$ for every $x\in L$.
A {\em sublattice} of $L$ is a submodule $L'$ of $L$, such that the induced
bilinear pairing on $L'$ is nondegenerate. A sublattice $L'$ of $L$ is
called {\em primitive} if $L/L'$ is torsion-free. The positive or
negative definiteness or signature of a lattice is defined to be that 
of the vector space $L_{\Q}$, together with the induced pairing.
\end{definition}

\begin{definition}\label{matLn}
For a lattice $L$ with pairing $\langle \underline{\,\,\,\,}\, ,
\underline{\,\,\,\,} 
\rangle$, we denote by $L(n)$ the lattice with the same underlying
module as $L$ and the pairing $n\cdot \langle \underline{\,\,\,\,}\, ,
\underline{\,\,\,\,} \rangle$. 
\end{definition}

\begin{definition}
The {\em Gram matrix} of a lattice $L$ with respect to a given basis
$x=(x_1,\ldots, x_n)$ is $I_x = (\langle x_i,x_j \rangle)_{i,j}$.
The {\em discriminant} of $L$ is defined by
$\disc L = \det I_x$ for any basis $x$ of
$L$. A lattice $L$ is called {\em unimodular} if it is integral and
$\disc L =\pm 1$.
\end{definition}

\begin{lemma}\label{matdetsub}
Let $L'$ be a sublattice of finite index in a lattice $L$. Then 
we have $\disc L' = [L:L']^2 \disc L$. 
\end{lemma}
\begin{proof}
This is a well known fact, see also \myref{\shioda}, section 6.
\end{proof}

\noindent
The definition of elliptic surface and the results in this section can
all be found in \myref{\shioda}. For a more detailed summary of these
results and constructions of elliptic surfaces, see also 
\myref{\luijkheron}, sections 3 and 4. Throughout this paper we will 
say that a variety $V$ over a field $k$ is smooth if the map $V
\rightarrow \Spec k$ is smooth.

\begin{definition}
Let $C$ be a smooth, irreducible, projective curve over an
algebraically closed field $k$. An {\em elliptic surface} over $C$ is
a smooth, irreducible, projective surface $S$, together with a
non-smooth, relatively minimal,
surjective morphism $f \colon S \ra C$, of which almost all fibers are
nonsingular curves of genus $1$, and a section $\O$ of $f$.
\end{definition}

\begin{remark}\label{matcastel}
By Castelnuovo's criterion (see \myref{\chin}, Thm. 3.1), the morphism
$f$ is relatively minimal if and only if no fiber contains an exceptional
divisor, i.e., a prime divisor $E$ with $E^2=-1$ and $H^1(E, \O_E)=0$.
By \myref{\hag}, Prop. III.9.7, any dominating morphism from an
integral variety to a regular curve is flat. 
Therefore, so is $f$ in the definition above. 
Also, $f$ is locally of finite presentation. Hence, by
\myref{\egatwo}, D\'ef. 6.8.1, the requirement of $f$ not being smooth
in the definition above is equivalent to the requirement that $f$ has
a singular fiber.
\end{remark}

\noindent
For the rest of this section,
let $S$ be an elliptic surface over a smooth, irreducible, projective curve
$C$ over an algebraically closed field $k$, fibered by $f\colon S \ra C$ 
with a section $\O$. Let $K=k(C)$ denote the function field of $C$ 
and let $\eta\colon \Spec K \rightarrow C$ be its generic point. 
Then the generic fiber 
$E=S \times_C \Spec K$ of $f$ is a smooth, projective, geometrically
integral curve over $K$ with 
genus $1$. Let $\xi$ denote the natural map $E\rightarrow S$. 

$$
\xymatrix{ E \ar[r]^{\xi} \ar[d] & S \ar[d]^f \\
           \Spec K \ar[r]_{\eta} &C}
$$

\begin{lemma}\label{matsections}
Both maps $\xi_*$ and $\eta^*$ in 
$$
E(K)=\Hom_K(\Spec K,E) \stackrel{\xi_*}{\longrightarrow}
\Hom_C(\Spec K,S) \stackrel{\eta^*}{\longleftarrow}
\Hom_C(C,S)=S(C)
$$
are bijective.
\end{lemma}
\begin{proof}
By the universal property of fibered products, we find that every morphism
$\sigma\colon \Spec K \rightarrow S$ with $f \circ \sigma = \eta$
comes from a unique section of the morphism $E \rightarrow \Spec
K$. Hence, the map $\xi_*$ is bijective. 
As $C$ is a smooth curve and $S$ is projective, any morphism from a
dense open subset of $C$ to $S$ extends uniquely to a morphism from
$C$, see \myref{\hag}, Prop. I.6.8. As $\Spec K$ is dense in $C$, the
map $\eta^*$ is bijective as well.
\end{proof}

\noindent
Whenever we implicitly identify the two sets $E(K)$
and $S(C)$, it
will be done using the bijection $\xi_*^{-1}\circ\eta^*$ of 
Lemma \ref{matsections}.
The section $\O$ of $f$ corresponds to a point on $E$, giving $E$ 
the structure of an elliptic curve. This endows $E(K)$ with a group
structure, which carries over to $S(C)$, see \myref{\silv}, Prop. III.3.4.  

Recall that for any proper scheme $Y$ over an algebraically closed field, 
the N\'eron-Severi group $\NS(Y)$ of $Y$ is the quotient of $\Pic Y$ by
the group $\Pic^0 Y$ consisting of all divisor classes algebraically 
equivalent to $0$, see \myref{\hag}, exc.~V.1.7, and \myref{\sgasix}, 
Exp. XIII, p. 644, 4.4. If $Y$ is proper, then $\NS(Y)$ is a
finitely generated, abelian group, see \myref{\hag}, exc.~V.1.7-8,
for surfaces, or \myref{\sgasix}, Exp. XIII, Thm. 5.1 in general.
Its rank $\rho=\dim \NS (Y) \otimes \Q$ is called the Picard number of $Y$.
Note that for the rest of this section $S$ is still an elliptic surface.

\begin{proposition}\label{matneronfree}
On $S$ algebraic equivalence coincides with numerical equivalence. The group
$\NS(S)$ is free. The intersection pairing induces a symmetric
nondegenerate bilinear pairing on $\NS(S)$, making it into a lattice
of signature $(1,\rho-1)$. If $S$ is a K3 surface, then $\NS(S)$ is an
even lattice.
\end{proposition}
\begin{proof}
The first statement is proved by Shioda in \myref{\shioda}, Thm.~3.1.
It follows immediately that the
bilinear intersection pairing is nondegenerate on $\NS(S)$, see
\myref{\shioda}, Thm.~2.1 or \myref{\hag}, example V.1.9.1. The
signature is given by the Hodge Index Theorem (\myref{\hag},
Thm.~V.1.9). If $S$ is a K3 surface, then its canonical sheaf is
trivial and the adjunction formula (\myref{\hag}, Prop. V.1.5) 
reduces to $D^2 = 2g(D)-2$ for any irreducible curve $D$ on $S$ with
genus $g(D)$. As the irreducible divisors generate $\NS(S)$, the
lattice $\NS(S)$ is even. 
%
% Actually, it is enough for the Euler characteristic of S to be even!
%
%
\end{proof}

\begin{lemma}\label{matsamepizero}
The induced map $f^* \colon \Pic^0 C \ra \Pic^0 S$ is an isomorphism. 
\end{lemma}
\begin{proof}
See \myref{\shioda}, Thm.~4.1. 
\end{proof}

\noindent
For every point $P \in E(K)$, let $(P)$ denote the prime divisor on $S$ 
that is the image of the section $C \ra S$ corresponding to $P$ by Lemma 
\ref{matsections}. Let $T \subset \NS(S)$ be generated
by the classes of the divisor $(\O)$ and the irreducible components of 
the singular fibers of $f$. For every $v \in C$, let $m_v$ 
denote the number of irreducible components of the fiber of $f$ at
$v$. Finally, let $r$ denote the rank of the Mordell-Weil group $E(K)$.

\begin{lemma}\label{matrangT}
The module $T$ is a sublattice of $\NS(S)$ of rank $\rk T = 2+ \sum_v
(m_v-1)$ and signature $(1,\rk T -1)$.
\end{lemma}
\begin{proof}
See \myref{\shioda}, Prop. 2.3.
\end{proof}

\begin{proposition}\label{matdoel}
There is a natural homomorphism $\varphi \colon \NS(S) 
\ra E(K)$ with kernel $T$.  It is surjective and maps $(P)$ to $P$. 
We have $\rho = \rk \NS(S) = r+2+ \sum_v (m_v-1)$.
\end{proposition}
\begin{proof}
The map $\varphi$ is defined in section 5 of \myref{\shioda}.
For surjectivity, see \myref{\shioda}, Lemma 5.1 and 5.2. The
fact that $T$ is the kernel is \myref{\shioda}, Thm. 1.3. The last
equality follows from Lemma \ref{matrangT} and the fact that the 
alternating sum of the ranks 
of finitely generated, abelian groups in an exact sequence equals $0$.
\end{proof}

\begin{corollary}\label{matEKmodTors}
There is a unique section $\psi$ of the homomorphism $\NS(S) \otimes
\Q \rightarrow E(K) \otimes \Q$ induced by $\varphi$ that maps $E(K)
\otimes \Q$ onto the orthogonal complement of $T \otimes \Q$ in
$\NS(S) \otimes \Q$. The homomorphism $\psi$ induces a symmetric
bilinear pairing on $E(K)$. The opposite of this pairing induces the
structure of a positive definite lattice on $E(K)/E(K)_{\tors}$.
\end{corollary}
\begin{proof}
See \myref{\shioda}, Thm. 8.4.
\end{proof}

\begin{remark}\label{matheightpair}
Shioda gives an
explicit formula for the pairing on $E(K)$, based on how the sections 
intersect the singular fibers and each other, see \myref{\shioda}, 
Thm. 8.6. 
\end{remark}
\bigskip

\mysection{\titlethree}\label{matsecthree}

\noindent
Let $G \subset \Aut X$ be the group of automorphisms of $X$ generated by
permutations of $x$, $y$ and $z$, by permutations of $a$, $b$, and $c$ and by
switching the sign of two of the coordinates $a$, $b$, and $c$. Then $G$ is
isomorphic to $(V_4 \semi S_3)\times S_3$ and has order $144$.
The surface $X$ has $12$ singular points, on which $G$ acts transitively.
They are all ordinary double points 
and their orbit under $G$ is represented by $[x:y:z:a:b:c]=[2:-1:-1:1:1:1]$.
Let $\pi \colon \bX \ra X$ be the blow-up of $X$ in these $12$ points.

Note that a K3 surface is a smooth, projective, geometrically
irreducible surface $S$, of which the canonical sheaf is trivial and
the irregularity $q=q(S)=\dim H^1(S, \O_S)$ equals $0$.

\begin{proposition}\label{matK3}
The surface $\bX$ is a smooth K3 surface. The exceptional curves 
above the $12$ singular points of $X$ are all isomorphic to $\P^1$ and have
self-intersection number $-2$. 
\end{proposition}
\begin{proof}
Ordinary double points are resolved after one blow-up, so $\bX$ is
smooth. The
exceptional curves $E_i$ are isomorphic to $\P^1$, see \myref{\hag},
exc. I.5.7. Their self-intersection number follows from \myref{\hag},
example V.2.11.4. 
Since $X$ is a complete intersection, it is geometrically connected and
$H^1(X,\O_X)=0$, so $q(X)=0$, see \myref{\hag}, exc. II.5.5. 
From its connectedness
it follows that $\bX$ is geometrically connected as well. As $\bX$ is
also smooth, it follows that $\bX$ is geometrically irreducible.

To compute the canonical
sheaf on $\bX$, note that on the nonsingular part $U=X^{\reg}$ of $X$ the
canonical sheaf is given by $\O_X(-5-1+3+2+1)|_U =
\O_{U}$, see \myref{\hag}, Prop. II.8.20 and exc. II.8.4. 
Hence, the canonical
sheaf on $\bX$ restricts to the structure sheaf outside the
exceptional curves. That implies that there are integers $a_i$ such
that $K = \sum_i a_i E_i$ is a canonical divisor. Recall that $E_i^2 =
-2$ and $E_i\cdot E_j =0$ for $i \neq j$. Applying the adjunction
formula $2g_C-2 = C \cdot (C+K)$ (see \myref{\hag}, Prop. V.1.5) to 
$C = E_i$, we find that $a_i=0$ for all $i$, whence $K=0$. 

It remains to show that $q(\bX)=q(X)$. 
It follows immediately from \myref{\artin}, Prop. 1,
that ordinary double points are rational singularities, i.e., we have 
$R^1 \pi_* \O_{\bX} = 0$. 
Also, as $X$ is integral, the sheaf $\pi_* \O_{\bX}$ is a 
sub-$\O_X$-algebra of the constant $\O_X$-algebra $K(X)$,
where $K(X)=K(\bX)$ is the function field of both $X$ and $\bX$. Since
$\pi$ is proper, $\pi_* \O_{\bX}$ is finitely generated as
$\O_X$-module. As $X$ is normal, i.e., $\O_X$ is integrally closed, we
get $\pi_* \O_{\bX} \isom \O_X$. Hence, the desired equality
$q(\bX)=q(X)$ follows from the following lemma, applied to $f=\pi$ and
$\sF = \O_{\bX}$.
\end{proof}

\begin{lemma}
Let $f \colon W \rightarrow Z$ be a continuous map of topological
spaces. Let $\sF$ be a sheaf of groups on $W$ and assume that 
$R^i f_*(\sF)=0$ for all $i=1, \ldots, t$. Then for all $i=0,1,\ldots,
t$, there are isomorphisms
$$
H^i(W, \sF) \isom H^i(Z, f_* \sF).
$$
\end{lemma}
\begin{proof}
This follows from the Leray spectral sequence. For a more elementary
proof, choose an injective resolution
$$
0 \ra \sF \ra I_0 \ra I_1 \ra I_2 \ra \cdots
$$
of $\sF$. Because $R^i f_*(\sF)=0$ for $i=1,\ldots, t$, we conclude
that the sequence
\begin{equation}\label{matflasque}
0 \ra f_* \sF\ra f_*I_0 \ra f_*I_1 \ra f_*I_2 \ra \cdots \ra f_* I_{t+1}
\end{equation}
is exact as well. As injective sheaves are flasque (see \myref{\hag},
Lemma III.2.4) and $f_*$ maps flasque $W$-sheaves to flasque
$Z$-sheaves, the exact sequence (\ref{matflasque}) can be extended to a
flasque resolution of $f_* \sF$. By \myref{\hag}, Rem. III.2.5.1, we
can use that flasque resolution to compute the cohomology groups 
$H^i(Z, f_* \sF)$. Taking global sections we get the complex
\begin{equation}\label{matglobalsectionscomplex}
0 \ra \Gamma(Z,f_*I_0) \ra \Gamma(Z,f_*I_1) \ra
\Gamma(Z,f_*I_2)  \ra \cdots \ra \Gamma(Z,f_*I_{t+1})\ra \ldots
\end{equation}
As $\Gamma(Z, f_*I_n) \isom \Gamma(W, I_n)$ for all $n$, we find that for
$i=0,1,\ldots, t$, the $i$-th cohomology of (\ref{matglobalsectionscomplex}) is
isomorphic to both $H^i(Z, f_*\sF)$ and $H^i(W,\sF)$.
\end{proof}

\noindent
We will now give $\bXbar$ the structure of an elliptic surface over $\P^1$.
Let $f \colon \bX \ra \P^1$ be the composition of $\pi$ with the morphism 
$f'\colon X \ra \P^1, [x:y:z:a:b:c] \mapsto [x:a] = [ 2bc : yz]$.
One easily checks that $f'$, and hence $f$, is well-defined everywhere.

If $a=0$, then clearly $M_{a,b,c}$ in (\ref{matmatrix}) has eigenvalue
$0$. Geometrically, this reflects the fact that the hyperplane $a=0$
intersects $X$ in three conics, one in each of the hyperplanes
given by $xyz=0$. Hence, each of the hyperplanes $H_t$ given by $x=ta$
in the family of hyperplanes through the space $x=a=0$ contains the
conic given by $a=x=0$ on $X$. The fibers of $f$ consist of the
inverse image under $\pi$ of the other components in the intersection
of $X$ with the family of hyperplanes $H_t$. The fiber above $[t:1]$ 
is therefore given by the intersection of the two quadrics 
\begin{equation}\label{matmatfibers}
xy+yz+zx = -a^2-b^2-c^2 \qquad \mbox{and} \qquad tyz = 2bc
\end{equation}
within the intersection of two
hyperplanes 
\begin{equation}\label{matinp3}
x+y+z=x-ta=0,
\end{equation}
which is isomorphic to $\P^3$. The conic $C$ given by $a+b=c-y=0$ on $X$
maps under $f'$ isomorphically to $\P^1$. The strict transform of $C$
on $\bX$ gives a section of $f$ that we will denote by $\O$.

\begin{proposition}
The morphism $f$ and its section $\O$ give $\bX_{\C}$ the structure of
an elliptic surface over $\P^1_\C$.
\end{proposition}
\begin{proof}
Since $\bX$ is a K3 surface, it is minimal. Indeed, by the adjunction
formula any smooth curve $C$ of genus $0$ on $\bX$ would have
self-intersection $C^2=-2$, while an exceptional curve that can be
blown down has self-intersection $-1$, see \myref{\hag}, Prop. V.3.1. 
Hence, $f$ is a relatively minimal fibration by Remark \ref{matcastel}. 
The $12$ exceptional curves give extra components in the
fibers above $t = \pm 1, \pm 2$, so $f$ is not smooth. 
From the description (\ref{matmatfibers})
above, an easy computation shows that the fibers above $t \neq 0, \pm
1, \pm 2, \infty$ are nonsingular. They are isomorphic to the complete 
intersection of two quadrics in $\P^3$, so by \myref{\hag},
exc. II.8.4g, almost all fibers have genus $1$.
\end{proof}

\noindent
Let $K\isom \Q(t)$ denote the function field of $\P^1_{\Q}$ and let
$E/K$ be the 
generic fiber of $f$. It can be given by the same equations (\ref{matmatfibers}) 
and (\ref{matinp3}). To put $E$ in Weierstrass form, set 
$\lambda=(t^2-4)\nu+3t$ and
$\mu=t(t^2-4)(z-y)(t\nu^2-2\nu+t)/x$, where $\nu=(x-c)/(a+b)$.
Then the change of variables
$$
\eqalign{
u &= \big( \mu+(\lambda^2+t(t^2-1)(t+8))\big)/2,\cr
v &= \big( 
\mu\lambda+\lambda^3+(t^2-1)(t^2-8)\lambda-8t(t^2-1)^2\big)/2\cr
}
$$
shows that $E/K$ is isomorphic to the elliptic curve over $K$ given by 
$$
v^2 = u\big(u-8t(t^2-1)\big)\big(u-(t^2-1)(t+2)^2\big).
$$
It has discriminant $\Delta = 2^{10} t^2(t^2-1)^6(t^2-4)^4$ and $j$-invariant
$$
j=\frac{4(t^4+56t^2+16)^3}{t^2(t^2-4)^4}.
$$
%   X = u,    Y = v,   Z = 1
%
%a=64*t^2*(t+2)^2*(-1+t)^3*(1+t)^4*Z^3-16*(-1+t)^2*(t+2)^2*(1+t)^2*(X*t-Y)*Z^2
%                       -8*X*(-1+t)*(1+t)*(-X*t+3*X*t^2+2*Y)*Z+2*t*X^3, 
%b=-64*t*(-2+t^2)*(t+2)^2*(-1+t)^3*(1+t)^4*Z^3+
%         16*t*(t+2)*(-1+t)^2*(1+t)^2*(-8*X-6*X*t+4*X*t^2+X*t^3-Y*t)*Z^2
%        -4*t*X*(-1+t)*(1+t)*(-10*X-6*X*t+3*X*t^2-Y*t)*Z-2*t*X^3, 
%c=-64*t^2*(t+2)^2*(-1+t)^3*(1+t)^4*Z^3
%       -16*t*(t+2)*(-1+t)^2*(1+t)^2*(-4*X-6*X*t-X*t^2+Y*t+2*X*t^3)*Z^2+
%       4*X*(-1+t)*(1+t)*(-4*X-12*X*t-4*X*t^2+6*X*t^3+2*Y*t+X*t^4)*Z-2*(-2+t^2)*X^3, 
%x=64*t^3*(t+2)^2*(-1+t)^3*(1+t)^4*Z^3-16*t*(-1+t)^2*(t+2)^2*(1+t)^2*(X*t-Y)*Z^2
%       -8*t*X*(-1+t)*(1+t)*(-X*t+3*X*t^2+2*Y)*Z+2*t^2*X^3, 
%y=-128*t*(t+2)^2*(-1+t)^3*(1+t)^4*Z^3
%      -32*(t+2)*(-1+t)^2*(1+t)^2*(-2*X-5*X*t-2*X*t^2+Y*t+2*X*t^3)*Z^2+
%      8*X*(-1+t)*(1+t)*(-6*X-6*X*t+X*t^2+3*X*t^3+Y*t)*Z-2*(-2+t^2)*X^3, 
%z=-64*t*(-2+t^2)*(t+2)^2*(-1+t)^3*(1+t)^4*Z^3+
%     16*(t+2)*(-1+t)^2*(1+t)^2*(-4*X-10*X*t-2*X*t^2+5*X*t^3-Y*t^2)*Z^2
%     -8*X*(-1+t)*(1+t)*(-6*X-6*X*t+2*X*t^2-Y*t)*Z-4*X^3

\begin{lemma}\label{matmatsingfibs}
The singular fibers of $f$ are at $t=0,\pm 1, \pm 2$ and at $t=\infty$. 
They are described in the following table, where $m_t$ (resp. $m_t^{(1)}$) 
is the number of irreducible components (resp. irreducible components
of multiplicity $1$). \newline \smallskip

\begin{center}
\begin{tabular}{|l|ccc|}
\hline
$t$        & type & $m_t$ & $m_t^{(1)}$ \cr
\hline
$0,\infty$ & $I_2$   & $2$ & $2$ \cr
$\pm 1$    & $I_0^*$ & $5$ & $4$ \cr
$\pm 2$    & $I_4$   & $4$ & $4$ \cr
\hline
\end{tabular}\smallskip
\end{center} 
\end{lemma}
\begin{proof}
This is a straightforward computation.
Since we have a Weierstrass form, it also follows easily from 
Tate's algorithm, see \myref{\tate} and \myref{\silvtwo}, IV.9.
\end{proof}

\noindent
Applying the automorphisms
$(b,c)\mapsto(-c,-b)$ and $(b,c)\mapsto(-b,-c)$ and 
$(b,c,y,z)\mapsto (c,b,z,y)$ to the curve $\O$, we get three more 
sections, which we will denote by $P$, $T_1$ and $T_2$ respectively.
By Lemma \ref{matsections}, these sections correspond with points on the
generic fiber $E/K$.
The Weierstrass coordinates $(u,v)$ of these points are given by 
\begin{equation}
\eqalign{
T_1 & = \big((t^2-1)(t+2)^2, \,0 \big),\cr
T_2 & = \big(0,\,0\big),\cr
P &= \big(2t^3(t+1), \,2t^2(t+1)^2(t-2)^2\big), \cr
}
\end{equation}

\noindent
We immediately notice that the $T_i$ are $2$-torsion points. 

\begin{proposition}\label{Pinford}
The section $P$ has infinite order in the group $S(C) \isom E(K)$.
\end{proposition}
\begin{proof}
Note that $S(C)$ and $E(K)$ are isomorphic by the identification of
Lemma \ref{matsections}.
By Corollary \ref{matEKmodTors} there is a bilinear 
pairing on $E(K)$ that induces a nondegenerate pairing 
on $E(K)/E(K)_{\tors}$. As mentioned in Remark \ref{matheightpair},
Shioda gives an explicit formula for this pairing, 
see \myref{\shioda}, Thm. 8.6.
We find that $\langle P,P\rangle = \frac{3}{2} \neq 0$, so $P$ is 
not torsion.
\end{proof}

\noindent
The main theorem now follows immediately. \newline \smallskip

\begin{proofof}{\bf Theorem \ref{matmain}.}
By Proposition \ref{Pinford}
the multiples of $P$ give infinitely many rational curves on $\bX$, so
the rational points on $\bX$ are dense. As $\pi$ is dominant, 
the rational points on $X$ are dense as well.
\end{proofof}

\noindent
The multiples of
$P$ yield infinitely many parametrizations of integral, symmetric 
$3\times 3$ matrices with zeros on the diagonal and integral eigenvalues. 
The section $2P$, for example, is a curve of
degree $8$ on $X$ which can be parametrized by 
$$
\eqalign{
a&=t(t^6-8t^4+20t^2-12), \cr
b&=-t(t^6-4t^4+4), \cr
c&=(t^2-2)(t^6-6t^4+8t^2-4), \cr
}
$$
and suitable polynomials for $x$, $y$, and $z$. The parametrization
(\ref{matlowdegparam}) does not come from a section of $f$. We will see
in Section 6 where it does come from.
\bigskip

\mysection{\titlefour}\label{matsecfour}

\noindent
As mentioned in the introduction, the geometry and the arithmetic 
of K3 surfaces are closely related. In the following sections we will 
further analyze the geometry of $\bX$. 
Set $L=\C(t) \supset \Q(t) = K$.
In this section we will find explicit generators for the
Mordell-Weil group $E(L)$ and for the 
N\'eron-Severi group of $\bXbar = \bX_{\C}$. 
This will be used in Sections 5 and 6. 

For any complex surface $Z$, the N\'eron-Severi group of $Z$ 
can be embedded in 
$H^{1,1}(Z)=H^1(Z,\Omega_Z^1)$, see \myref{\BPV}, p. 120. 
If $Z$ is a complex K3 surface, we have $\dim H^{1,1}(Z)=20$, see 
\myref{\BPV}, Prop. VIII.3.3. Hence we find that the Picard number
$\rho(Z) = \rk \NS(Z)$ is at most $20$. If $\rho(Z)$ is equal to $20$
we say that $Z$ is a singular K3 surface.

\begin{proposition}
The Picard group $\Pic \bXbar$ is isomorphic to $\NS(\bXbar)$ and it is a
finitely generated, free abelian group.
\end{proposition}
\begin{proof}
As $\bXbar$ has the structure of an elliptic surface over $\P^1$ and
$\Pic^0 \P^1 =0$, the isomorphism follows from Lemma \ref{matsamepizero}.
The last statement follows from Proposition \ref{matneronfree}.
\end{proof}

\noindent
Two of the irreducible components of the singular fibers of $f \colon
\bX \ra \P^1$ above $t =
\pm 2$ are defined over $\Q(\sqrt{3})$. They are all in the same orbit
under $G$. In that same orbit we also find a section, given by 
$z=2b$ and $2(c-a)=\sqrt{3}(y-x)$. We will denote it by $Q$. Its
Weierstrass coordinates are given by 
$$
Q = \big(2t(t+1)(t+2),\, 2\sqrt{3}t(t^2-4)(t+1)^2 \big).
$$
It follows immediately that the Galois conjugate of $Q$ 
under the automorphism that sends $\sqrt{3}$ to $-\sqrt{3}$ is equal to $-Q$.

\begin{proposition}\label{matmatgroup}
The surface $\bXbar$ is a singular K3 surface. The Mordell-Weil group
$E(L)$ is 
isomorphic to $\Z^2\times (\Z/2\Z)^2$ and generated by $P$, $Q$, $T_1$ and
$T_2$. The Mordell-Weil group $E(K)$ is
isomorphic to $\Z \times (\Z/2\Z)^2$ and generated by $P$, $T_1$ and $T_2$.
\end{proposition}
\begin{proof}
From Shioda's explicit formula for the pairing on $E(K)$
(see Remark \ref{matheightpair}), we find that 
$\langle P,P \rangle = \frac{3}{2}$ and 
$\langle Q,Q \rangle = \frac{1}{2}$ and $\langle P,Q \rangle = 0$.
Hence, $P$ and $Q$ are linearly independent and the Mordell-Weil rank 
$r=\rk E(L)$ is at least $2$.

By Lemma \ref{matmatsingfibs} and \ref{matrangT}, the lattice  
generated by the vertical fibers and $\O$ has rank $18$. From Proposition 
\ref{matdoel} it follows that the rank $\rho$ of $\NS(\bXbar)=\Pic(\bXbar)$ is 
at least $18+2=20$. As $\bXbar$ is a K3 surface (see Proposition \ref{matK3})
and $20$ is the maximal Picard number for K3 surfaces in
characteristic $0$, we conclude that $\bXbar$ is a singular K3 surface. 
Using Proposition \ref{matdoel} again, we find that the 
Mordell-Weil rank of $E(L)$ equals $2$. 
Since $E$ has additive reduction at $t=\pm 1$, the order of 
the torsion group $E(L)_{\tors}$ is at most $4$, see \myref{\silvtwo}, Remark
IV.9.2.2. Hence we have $E(L)_{\tors} = \langle T_1,T_2 \rangle$. 

From Shioda's explicit formula for the height pairing it follows 
that with singular fibers only of type $I_2$, $I_4$ and
$I_0^*$, the pairing takes values in $\frac{1}{4}\Z$.
Hence, the lattice $\Lambda = \big(E(L)/E(L)_\tors\big)(4)$ is
integral, see Definition \ref{matLn}. In $\Lambda$ we have $\langle
P,P\rangle = 6$ and $\langle Q,Q \rangle = 2$ and 
$\langle P,Q \rangle = 0$. Hence, by Lemma \ref{matdetsub} 
the sublattice $\Lambda'$ of $\Lambda$ 
generated by $P$ and $Q$ has discriminant $\disc \Lambda'=12=n^2 \disc
\Lambda$, with $n=[\Lambda:\Lambda']$.
Therefore, $n$ divides $2$. Suppose $n=2$. Then there is an $R \in \Lambda
\setminus \Lambda'$ with $2R=aP+bQ$. By adding multiples of $P$ and
$Q$ to $R$, we may assume $a,b \in \{0,1\}$. In $\Lambda$ we get 
$4| \langle 2R,2R \rangle = 6a^2+2b^2$. Hence,
we find $a=b=1$, so $2R=P+Q+T$ for some torsion element $T \in E(L)[2]$.
Since all the $2$-torsion of $E(L)$ is rational over $L$, it is 
easy to check whether an element of $E(L)$ is in
$2E(L)$. If $e$ is
the Weierstrass $u$-coordinate of one of the $2$-torsion points, 
then there is a homomorphism 
$$
E(L)/2E(L) \ra L^*/{L^*}^2,
$$
given by $S \mapsto u(S)-e$, where $u(S)$ denotes the Weierstrass
$u$-coordinate of the point $S$, see \myref{\silv}, \parag{} X.1.
We can use $e=0$ and find that for none of the four torsion points 
$T \in E(L)[2]$ the value $u(P+Q+T)$ is a square in $L$. 
Hence, we get $n=1$ and $E(L)$ is generated by $P,Q,T_1,$ and $T_2$.

Suppose $aP+bQ+\varepsilon_1T_1+\varepsilon_2T_2$ is contained in $E(\Q(t))$
for some integers $a,b,\varepsilon_i$. Then also $bQ\in E(\Q(t))$.
As the Galois automorphism $\sqrt{3} \mapsto -\sqrt{3}$ sends $Q$ to
$-Q$, we find that $bQ= -bQ$.  But $Q$ has infinite order, so $b=0$. 
Thus, we have $E(\Q(t)) = \langle P,T_1,T_2\rangle$.
\end{proof}

\noindent
To work with explicit generators of the N\'eron-Severi group of
$\bXbar$, we will
name some of the irreducible divisors that we encountered so far
as in the table below. The exceptional curves are given by the
point on $\Xbar=X_{\C}$ that they lie above. 
Other components of singular fibers
are given by their equations on $\Xbar$. Sections are given by their
equations and the name they already have.

$$
\begin{array}{llll}
D_1  & x=-2a, b+c=\frac{\sqrt{3}}{2}(y-z) & D_{11}      &
[-1:-1:2:-1:-1:1] \cr 
D_2  & [2:-1:-1:-1:1:-1] & D_{12}      & (T_1) \colon a-b=c+y=0  \cr
D_3  & (\O) \colon a+b=c-y=0 & D_{13}      & [2:-1:-1:1:1:1] \cr
D_4  & [-1:-1:2:1:-1:-1] & D_{14}      & x=2a, 2(b-c)=\sqrt{3}(y-z) \cr
D_5  & a=-x, b=c & D_{15}      & (Q) \colon z=2b, 
                                          c-a=\frac{\sqrt{3}}{2}(y-x) \cr
D_6  & [-1:2:-1:1:-1:-1] & D_{16}      & x=2a, 2(b-c)= \sqrt{3}(z-y) \cr
D_7  & (T_2) \colon a+c=b-z=0 & D_{17}      & x=b=0   \cr
D_8  & [-1:2:-1:1:1:1] & D_{18}      & a=y=0   \cr
D_9  & [-1:2:-1:-1:1:-1] & D_{19}      & (P) \colon a-c=b+y=0    \cr
D_{10}&a=x, b=-c & D_{20}      & F \text{ (whole fiber)} \cr
\end{array}
$$
\smallskip

\begin{proposition}\label{matgrammprop}
The sequence $\{D_1, D_2, \ldots, D_{20}\}$ forms an ordered basis for
the N\'eron-Severi lattice $\NS(\bXbar)$. With respect to this basis 
the Gram matrix of inner products is given by 
{\small
$$
\left(
\begin{array}{*{18}{r@{\,}}r@{\quad\,}r}
 -2& 1& 0& 0& 0& 0& 0& 0& 0& 0& 0& 0& 0& 0& 0& 0& 0& 0& 0& 0\\
  1&-2& 1& 0& 0& 0& 0& 0& 0& 0& 0& 0& 0& 0& 0& 0& 0& 0& 1& 0\\
  0& 1&-2& 1& 0& 0& 0& 0& 0& 0& 0& 0& 0& 0& 0& 0& 1& 0& 0& 1\\
  0& 0& 1&-2& 1& 0& 0& 0& 0& 0& 0& 0& 0& 0& 0& 0& 0& 0& 0& 0\\
  0& 0& 0& 1&-2& 1& 0& 1& 0& 0& 0& 0& 0& 0& 0& 0& 0& 0& 0& 0\\
  0& 0& 0& 0& 1&-2& 1& 0& 0& 0& 0& 0& 0& 0& 0& 0& 0& 0& 0& 0\\
  0& 0& 0& 0& 0& 1&-2& 0& 0& 0& 0& 0& 0& 0& 0& 0& 0& 1& 1& 1\\
  0& 0& 0& 0& 1& 0& 0&-2& 0& 0& 0& 0& 0& 0& 0& 0& 0& 0& 0& 0\\
  0& 0& 0& 0& 0& 0& 0& 0&-2& 1& 0& 0& 0& 0& 0& 0& 0& 0& 0& 0\\
  0& 0& 0& 0& 0& 0& 0& 0& 1&-2& 1& 0& 0& 0& 0& 0& 0& 0& 0& 0\\
  0& 0& 0& 0& 0& 0& 0& 0& 0& 1&-2& 1& 0& 0& 0& 0& 0& 0& 0& 0\\
  0& 0& 0& 0& 0& 0& 0& 0& 0& 0& 1&-2& 1& 0& 0& 0& 1& 0& 0& 1\\
  0& 0& 0& 0& 0& 0& 0& 0& 0& 0& 0& 1&-2& 1& 0& 1& 0& 0& 1& 0\\
  0& 0& 0& 0& 0& 0& 0& 0& 0& 0& 0& 0& 1&-2& 1& 0& 0& 0& 0& 0\\
  0& 0& 0& 0& 0& 0& 0& 0& 0& 0& 0& 0& 0& 1&-2& 0& 0& 1& 0& 1\\
  0& 0& 0& 0& 0& 0& 0& 0& 0& 0& 0& 0& 1& 0& 0&-2& 0& 0& 0& 0\\
  0& 0& 1& 0& 0& 0& 0& 0& 0& 0& 0& 1& 0& 0& 0& 0&-2& 0& 0& 0\\
  0& 0& 0& 0& 0& 0& 1& 0& 0& 0& 0& 0& 0& 0& 1& 0& 0&-2& 0& 0\\
  0& 1& 0& 0& 0& 0& 1& 0& 0& 0& 0& 0& 1& 0& 0& 0& 0& 0&-2& 1\\
  0& 0& 1& 0& 0& 0& 1& 0& 0& 0& 0& 1& 0& 0& 1& 0& 0& 0& 1& 0\\
\end{array}
\right).
$$
}
\end{proposition}
\begin{proof}
By Proposition \ref{matdoel} the N\'eron-Severi group 
$\NS(\bXbar)$ is
generated by $(\O)$, all irreducible components of the singular fibers, 
and any set of generators of the Mordell-Weil group $E(L)$. Thus, from 
Lemma \ref{matmatsingfibs} and Proposition \ref{matmatgroup} we can find a
set of generators for $\NS(\bXbar)$. 
Using a computer algebra package or even by hand, one checks that 
$\{D_1,\ldots, D_{20}\}$ generates the same lattice. A big part of the
Gram matrix is easy to compute, as we know how all fibral divisors
intersect each other. Also, every section intersects each fiber in
exactly one irreducible component, with multiplicity $1$. The sections
are rational curves, so by the adjunction formula they have
self-intersection $-2$. That leaves $\binom{5}{2}$
more unknown intersection numbers among the sections. By applying
appropriate automorphisms from $G\subset \Aut X$, we find that they 
are equal to intersection numbers that are already known by the above.
\end{proof}

\begin{remark}\label{matwhatshyp}
By Proposition \ref{matgrammprop} the hyperplane section $H$ is
numerically equivalent with a linear combination of the $D_i$. 
This linear combination is uniquely determined by the 
intersection numbers $H \cdot D_i$ for $i=1,\ldots, 20$ and turns out
to be some uninformative linear combination with many nonzero
coefficients. 
The reason for choosing the $D_i$ and their order in this manner is
that $D_1,\ldots, D_8$ and $D_9, \ldots, D_{16}$ generate
two orthogonal sublattices, both isomorphic to $E_8(-1)$.
In fact, we have the following proposition, which will be used in
Section 5.
\end{remark}

\begin{proposition}\label{matstructNS}
The N\'eron-Severi lattice $\NS(\bXbar)$ has discriminant $-48$. It
is isomorphic to the orthogonal direct sum
$$
E_8(-1) \oplus E_8(-1) \oplus \Z(-2) \oplus \Z(-24) \oplus U,
$$
where $U$ is the unimodular lattice with Gram matrix
$$
\left(
\begin{array}{cc}
0& 1 \\
1& 0 \\
\end{array}
\right)
$$
\end{proposition}
\begin{proof}
The discriminant of $\NS(\bXbar)$ is the determinant of the Gram
matrix, which equals $-48$. With respect to the basis $D_1, \ldots, D_{20}$, 
let $C_1, \ldots, C_4$ be defined by
$$
\eqalign{
C_1 = (& 0,0,0,-1,-2,-2,-2,-1,1,2,3,4,4,2,0,2,1,-2,0,0)\cr 
C_2 = (&6,12,26,29,32,19,6,16,9,18,27,36,34,23,12,17,7,-3,-8,4)\cr 
C_3 = (&1,2,4,4,4,2,0,2,2,4,6,8,8,5,2,4,2,-1,-1,0)\cr 
C_4 = (&1, 2, 4, 5, 6, 4, 2, 3, 1, 2, 3, 4, 4, 3, 2, 2, 0, 0, -1, 1)\cr
}
$$
and let $L_1,\ldots,L_5$ be the lattices generated by $(D_1,\ldots,D_8)$, 
$(D_9,\ldots,D_{16})$, $(C_1)$, $(C_2)$, and $(C_3,C_4)$
respectively. Then one easily checks that $L_1,\ldots,L_5$ are
isomorphic to $E_8(-1)$, $E_8(-1)$, $\Z(-2)$, $\Z(-24)$, and $U$
respectively. They are 
orthogonal to each other, and the orthogonal direct sum 
$L=L_1 \oplus \cdots \oplus L_5$ has discriminant $-48$ and rank
$20$. By Lemma \ref{matdetsub} we find that the index $[\NS(\bXbar):L]$
equals $1$, so $\NS(\bXbar)=L$.
\end{proof}

\bigskip

\mysection{\titlefive}\label{matsecfive}

\noindent
If $A$ is an abelian surface, then the involution $\iota = [-1]$ has
$16$ fixed points. The quotient $A/\langle \iota \rangle$ therefore
has $16$ ordinary double points. A minimal resolution of such a quotient 
is called a Kummer surface. All Kummer surfaces are K3 surfaces. 
Because of their rich geometric structure, their arithmetic can be
analyzed and described more easily. Every complex singular surface is
either a Kummer surface or a double cover of a Kummer surface, see
\myref{\SI}, Thm. 4 and its proof. 
It is therefore natural to ask whether our complex
singular K3 surface $\bXbar$ has the rich structure of a Kummer
surface. In Corollary \ref{matnonkummer} we will see that this is not the
case. 

Shioda and Inose have classified complex singular K3 surfaces by 
showing that the set of their isomorphism classes is
in bijection with the set of equivalence classes of positive definite 
even integral binary quadratic forms modulo the action of $\SL_2(\Z)$, see
\myref{\SI}. 
A singular K3 surface $S$ corresponds with the binary
quadratic form given by the intersection product on the oriented
lattice $T_S= \NS(S)^{\perp}$ of transcendental cycles on $S$.
Here the orthogonal complement is taken in the unimodular 
lattice $H^2(S,\Z)$ of signature $(3,19)$ (see \myref{\BPV}, Prop. VIII.3.2).
To find out which quadratic form the surface $\bXbar$ corresponds to, we
will use discriminant forms as defined by Nikulin \myref{\nik},
\parag{} 1.3.

\begin{definition}
Let $A$ be a finite abelian group.
A {\em finite symmetric bilinear form} on $A$ is a symmetric bilinear map 
$b\colon A \times A \rightarrow \Q/\Z$.

A {\em finite quadratic form} on $A$ is a map $q \colon A \rightarrow \Q/2\Z$,
such that for all $n \in \Z$ and $a \in A$ we have $q(na)=n^2q(a)$ and
such that the unique map $b \colon A \times A \rightarrow \Q/\Z$ determined
by $q(a+a')-q(a)-q(a') \equiv 2b(a,a') \mod 2\Z$ for all $a,a' \in A$
is a finite symmetric bilinear form on $A$.
The form $b$ is called {\em the bilinear form of $q$}. 
\end{definition}

\begin{definition}\label{defquadform}
Let $L$ be an integral lattice. We define the {\em dual lattice} $L^*$ by 
$$
\{x \in L_\Q \,\,|\,\, \langle x,y \rangle \in \Z \text{ for all } y \in L\}.
$$
\end{definition}

\begin{lemma}\label{matdetispreAL}
Let $L$ be an integral lattice. Then $|\!\disc L| = [L^*:L]$. 
\end{lemma}
\begin{proof}
There is an isomorphism $L^* \isom \Hom(L,\Z)$. If $x$ is a basis for
$L$, then the standard dual basis $x'$ of $\Hom(L_\Q,\Q)$ generates
$\Hom(L,\Z)$ as a $\Z$-module. Hence, for the Gram matrices $I_x$ and
$I_{x'}$ we find $I_{x'} = I_x^{-1}$. Thus, $\disc L^*=1/(\disc
L)$. By Lemma \ref{matdetsub} we have $\disc L = [L^*:L]^2 \disc L^*$, 
from which the equality follows. 
\end{proof}

\begin{lemma}\label{discform}\label{matdetisAL}
Let $L$ be an even lattice and set $A_L = L^*/L$. Then we have $\# A_L =
|\disc L|$ and the map
$$
q_L \colon A_L \rightarrow\Q/2\Z\colon x\mapsto\langle x,x\rangle+2\Z
$$
is a finite quadratic form on $A_L$.
\end{lemma}
\begin{proof}
The first statement is a reformulation of Lemma \ref{matdetispreAL}.
The map $q_L$ is well defined, as for $x \in L^*$ and $\lambda
\in L$, we have $\langle x+\lambda,x+\lambda \rangle - \langle x,x
\rangle = 2 \langle x, \lambda \rangle + \langle \lambda, \lambda
\rangle \in 2\Z$. The unique map $b\colon A_L \times A_L \ra \Q/\Z$ as
in Definition \ref{defquadform} is given by $(a,a') \mapsto \langle a,
a' \rangle +\Z$, which is clearly a finite symmetric bilinear
form. Thus, $q_L$ is a finite quadratic form.
\end{proof}
 
\begin{definition}\label{discquadform}
If $L$ is an even lattice, then the map $q_L$ as in
Lemma {\rm \ref{discform}} is called the {\em discriminant-quadratic form}
associated to $L$.
\end{definition}

\begin{lemma}\label{matquad}
Let $L$ be a primitive sublattice of an even unimodular lattice
$\Lambda$. Let $L^{\perp}$ denote the orthogonal complement of $L$ in
$\Lambda$. Then $q_L \isom -q_{L^{\perp}}$, i.e., there is an
isomorphism $A_L \rightarrow A_{L^{\perp}}$ making the following
diagram commutative.
$$
\xymatrix{
A_L \ar[r]^{\isom}\ar[d]_{q_L} & A_{L^{\perp}} \ar[d]^{q_{L^{\perp}}} \\
\Q/2\Z \ar[r]^{[-1]} & \Q/2\Z 
}
$$          
\end{lemma}
\begin{proof}
See \myref{\nik}, Prop. 1.6.1.
\end{proof}

\begin{lemma}\label{matwhereNS}
The embedding $\NS(\bXbar) \ra H^2(\bXbar,\Z)$ makes $\NS(\bXbar)$ into a
primitive sublattice of the even unimodular lattice $H^2(\bXbar,\Z)$.
We have $\disc T_{\bXbar}=48$.
\end{lemma}
\begin{proof}
For the fact that $H^2(\bXbar,\Z)$ is even and unimodular see \myref{BPV},
Prop. VIII.3.2. The image of the N\'eron-Severi group in $H^2(\bXbar,\Z)$
is equal to $H^{1,1}(\bXbar) \cap H^2(\bXbar,\Z)$, where the intersection is
taken in $H^2(\bXbar,\C)$, see \myref{\BPV}, p. 120. 
Hence, $\NS(\bXbar)$ is a primitive sublattice. From Lemma
\ref{matdetisAL} and \ref{matquad} we find 
$$
|\!\disc T_{\bXbar}| = |A_{T_{\bXbar}}| = |A_{\NS(\bXbar)}| = 
    |\!\disc \NS(\bXbar)| = 48.
$$
As $T_{\bXbar}$ is positive definite, we get $\disc T_{\bXbar}=48$. 
\end{proof}

\noindent
Up to the action of $\SL_2(\Z)$,
there are only four $2$-dimensional positive definite even lattices with 
discriminant $48$. The transcendental lattice $T_{\bXbar}$ is
equivalent to one of them. They are given by the Gram matrices
\begin{equation}\label{matfourmat}
\left(
\begin{array}{cc}
2&0 \cr
0&24 \cr
\end{array}
\right), \qquad
\left(
\begin{array}{cc}
4&0 \cr
0&12 \cr
\end{array}
\right), \qquad
\left(
\begin{array}{cc}
8&4 \cr
4&8 \cr
\end{array}
\right), \qquad
\left(
\begin{array}{cc}
6&0 \cr
0&8 \cr
\end{array}
\right). 
\end{equation}

\begin{proposition}
Under the correspondence of Shioda and Inose, the singular K3 surface $\bXbar$
corresponds to the matrix
$$
\left(
\begin{array}{cc}
2&0 \cr
0&24 \cr
\end{array}
\right).
$$
\end{proposition}
\begin{proof}
As $E_8(-1)$ and $U$ as in Proposition \ref{matstructNS} are unimodular,
it follows from Proposition \ref{matstructNS} and Lemma \ref{matdetisAL} that the
discriminant-quadratic form of $\NS(\bXbar)$ is isomorphic to that of
$\Z(-2)\oplus \Z(-24)$. By Lemma \ref{matquad} and \ref{matwhereNS} we find 
that the discriminant-quadratic form associated to $T_{\bXbar}$ is
isomorphic to that of $\Z(2)\oplus \Z(24)$, whence it takes on the
value $\frac{1}{24}+2\Z$. Of the four lattices described in 
(\ref{matfourmat}), the lattice $\Z(2)\oplus \Z(24)$ is the only one for
which that is true.
\end{proof}

\begin{corollary}\label{matnonkummer}
The surface $\bXbar$ is not a Kummer surface.
\end{corollary}
\begin{proof}
By \myref{\inose}, Thm. 0, a singular K3 surface $S$ 
is a Kummer surface if and 
only if its corresponding positive definite even integral binary
quadratic form 
is twice another such form, i.e., if $x^2 \equiv 0 \mod 4$ for all $x
\in T_S$. This is not true in our case.
\end{proof}
\bigskip

\mysection{\titlesix}\label{matsecsix}

\noindent
Note that so far we have seen $63$ rational curves of degree $2$ on $\Xbar$, 
namely those in the orbits under $G$ of 
\begin{equation}\label{matmatcurves}
\begin{array}{lll}
D_{10}\colon &x=a, &   b=-c, \cr
D_{16}\colon &x=2a, &  2(b-c) = \sqrt{3}(z-y), \cr
D_{17}\colon &x=0, &  b=0. \cr
\end{array}
\end{equation}
These orbits have sizes $18$, $36$, and $9$ respectively.
All of these curves correspond to infinitely many matrices 
that are either trivial or not defined over $\Q$.
To find more rational curves of low degree, we look at
fibrations of $\bXbar$ other than $f$. 
The conic $(\O)$ given by $a+b=c-y=0$ on $X$ determines 
a plane in the four-space in $\P^5$ given by $x+y+z=0$. The family 
of hyperplanes in this four-space
that contain that plane, cut out another family of elliptic curves on $\bX$. 
One singular fiber in this family is contained in the 
hyperplane section  $a+b=2(c-y)$ on $X$. It is the degree $4$
curve corresponding to the parametrization in (\ref{matlowdegparam}).
We will now see that this is the lowest degree of a parametrization of
nontrivial matrices defined over $\Q$.

Recall that $G \subset \Aut X$ is the group of automorphisms of 
$X$ generated by
permutations of $x$, $y$ and $z$, by permutations of $a$, $b$, and $c$ and by
switching the sign of two of the coordinates $a$, $b$, and $c$.

\begin{proposition}\label{matallcurves}
The union of the three orbits under the action of $G$ of the curves 
described in {\rm (\ref{matmatcurves})} consists of all $63$ curves 
on $\Xbar$ of degree smaller than $4$. 
\end{proposition}

Arguments similar to the ones used to prove Proposition \ref{matallcurves}
can be found in \myref{\bremner}, p.\ 302. To prove this final Proposition 
\ref{matallcurves} we will use the following lemma.

\begin{lemma}\label{matcurvesonK3}
Let $S$ be a minimal, nonsingular, algebraic K3 surface over $\C$.
Suppose $D$ is a divisor on $S$ with $D^2=-2$. \vspace{-2mm}
\begin{itemize}
\setlength{\itemsep}{-1mm}
\item[\rm (a)] If $D \cdot H$ is positive for some ample divisor $H$ on $S$,
then $D$ is linearly equivalent with an effective divisor.
\item[\rm (b)] If $D$ is effective and its corresponding closed subscheme
is reduced and simply connected, then the complete linear
system $|D|$ has dimension $0$. 
\end{itemize}
\end{lemma}
\begin{proof}
Since the canonical sheaf on $S$ is trivial and the Euler
characteristic $\chi$ of $\O_S$ equals $2$, the 
Riemann-Roch Theorem for surfaces (see \myref{\hag}, Thm V.1.6)
tells us that
$$
l(D)-s(D)+l(-D) = \frac{1}{2} D^2 + \chi = 1,
$$
where $l(D) = \dim H^0(S,\L(D)) = \dim |D|+1$ and $s(D) = \dim
H^1(S,\L(D))$ is the superabundance. For (a) it is enough to prove 
$l(D) \geq 1$. Because $s(D)$ is nonnegative, it suffices to show 
$l(-D)=0$. As we have $(-D)\cdot H < 0$, this follows from the fact 
that effective divisors have nonnegative intersection with ample divisors. 
For (b), $D$ is effective, so we also find $l(-D)=0$. 
In order to prove $l(D)=1$, it suffices to show that 
$s(D)=0$ or by symmetry, that $s(-D)=0$. Now $\L(-D)$ is equal to the
ideal sheaf $\I_Z$ of the closed subscheme $Z$ corresponding to $D$ and 
$H^1(S,\L(-D)) = H^1(S,\I_Z)$ fits in the exact sequence
$$
H^0(Z,\O_Z) \ra H^0(S,\O_S) \ra H^1(S,\I_Z) \ra H^1(Z,\O_Z).
$$
As $S$ and $Z$ are projective and connected, the first map is an 
isomorphism of one-dimensional vector spaces. Hence the map 
$H^1(S,\I_Z) \ra H^1(Z,\O_Z)$ is injective. By the Hodge decomposition
we know that $H^1(Z,\O_Z)$ is a direct 
summand of $H^1(Z,\C)$. Hence it is trivial, as $Z$ is simply
connected. Therefore, also $H^1(S,\I_Z)$ is trivial and $s(-D)=0$.
\end{proof}

\begin{proofof}{\bf Proposition \ref{matallcurves}.}
Let $C$ be a curve on $\Xbar$ of degree $d$ and arithmetic genus $g_a$
and let $C$
also denote its strict transform on $\bXbar$. Let its coordinates with
respect to the basis $\{D_1,\ldots, D_{20}\}$ of $\NS(\bXbar)$ 
be given by $m_1,\ldots,
m_{20}$. Let $H$ denote a hyperplane section. If $E$ is any of the $12$
exceptional curves on $\bXbar$, then we have $H \cdot E =0$. For any
curve $D$ on $\Xbar$ we have $H \cdot D = \deg D$. This
determines $H \cdot D_i$ for all $i=1,\ldots, 20$ (see Remark
\ref{matwhatshyp}), and we find 
\begin{equation}\label{matmatdegree}
\eqalign{
d = C\cdot H = 2\big(m_1&+m_3+m_5+m_7+m_{10}+m_{12}+m_{14}+\cr
&+m_{15}+m_{16}+m_{17}+m_{18}+m_{19}+2m_{20}\big).
}
\end{equation}
This implies that $d$ is even, say $d=2k$. Since we have $H^2=6$,
we can write the divisor class $[C] \in \NS(\overline{Y})$ as 
$[C] = \frac{d}{6} H + D = \frac{k}{3}H+D$ for some element $D\in
\frac{1}{6}\langle H \rangle^{\perp}$, where the orthogonal complement
is taken inside $\NS(\overline{Y})$. 
From the adjunction formula (see \myref{\hag}, Prop. V.1.5) we find
$C^2 = 2g_a-2$, so from $C^2=D^2+(\frac{kH}{3})^2$
we get $D^2 = 2g_a-2-\frac{2k^2}{3}$. 
By the Hodge Index Theorem (\myref{\hag}, Thm.~V.1.9) the lattice 
$\frac{1}{e} \langle H \rangle^{\perp}$ is negative definite for any $e>0$, 
so for fixed $k$ and $g_a$ there are only finitely 
many elements $D\in \frac{1}{6} \langle H \rangle^{\perp}$ with 
$D^2=2g_a-2-\frac{2k^2}{3}$. We will now make this more concrete. Set
%
%This implies that $d$ is even, say $d=2k$. We find 
%from the adjunction formula (see \myref{\hag}, Prop. V.1.5) that 
%$C^2 = 2g_a-2$. Furthermore, in the lattice $\NS(\bXbar)$, the element $C$ 
%is contained in the hyperplane given by $H\cdot D= d$. By the Hodge
%Index Theorem (\myref{\hag}, Thm.~V.1.9) this 
%coset of the orthogonal complement of $H$ is
%negative definite, so for fixed $d$ and $g_a$ there are only finitely 
%many elements $D$ with $D^2=2g_a-2$ and $H\cdot D=d$. 
%We will now make this more concrete. Set
%
{\footnotesize
$$
\eqalign{
v_{1} =& 2m_{2}+m_{5}+m_{7}+m_{10}+m_{12}+m_{14}+m_{15}+m_{16}+m_{17}+m_{18}+2m_{20}-k, \cr
v_{2} =&
4m_{3}-m_{4}+2m_{5}+2m_{7}+2m_{10}+2m_{12}+2m_{14}+2m_{15}+2m_{16}+m_{17}+\cr
       &+2m_{18}+2m_{19}+3m_{20}-2k, \cr
v_{3} =&
7m_{4}-2m_{5}+2m_{7}+2m_{10}+2m_{12}+2m_{14}+2m_{15}+2m_{16}+m_{17}+2m_{18}+\cr
       &+2m_{19}+3m_{20}-2k, \cr
v_{4} =&
33m_{5}-14m_{6}+9m_{7}-14m_{8}+9m_{10}+9m_{12}+9m_{14}+9m_{15}+9m_{16}+15m_{17}+\cr
       &+9m_{18}+16m_{19}+24m_{20}-9k, \cr
v_{5} =&
52m_{6}-24m_{7}-14m_{8}+9m_{10}+9m_{12}+9m_{14}+9m_{15}+9m_{16}+15m_{17}+9m_{18}+\cr
       &+16m_{19}+24m_{20}-9k, \cr
v_{6} =&
24m_{7}+m_{8}+4m_{10}+4m_{12}+4m_{14}+4m_{15}+4m_{16}+11m_{17}-9m_{18}-3m_{19}+\cr
       &+2m_{20}-4k, \cr
v_{7} =&
35m_{8}+8m_{10}+8m_{12}+8m_{14}+8m_{15}+8m_{16}+13m_{17}+9m_{18}+15m_{19}+\cr
       &+22m_{20}-8k, \cr
v_{8} =& 2m_{9}-m_{10}, \cr
v_{9} =&
211m_{10}-140m_{11}+m_{12}+m_{14}+m_{15}+m_{16}+41m_{17}+23m_{18}+50m_{19}+\cr
      &+64m_{20}-k, \cr
v_{10} =& 282m_{11}-210m_{12}+m_{14}+m_{15}+m_{16}+41m_{17}+23m_{18}+50m_{19}+64m_{20}-k, \cr
v_{11} =& 119m_{12}-94m_{13}+m_{14}+m_{15}+m_{16}-53m_{17}+23m_{18}+50m_{19}-30m_{20}-k, \cr
v_{12} =& 144m_{13}-118m_{14}+m_{15}-118m_{16}-53m_{17}+23m_{18}-69m_{19}-30m_{20}-k,\cr
v_{13} =& 86m_{14}-71m_{15}-58m_{16}-5m_{17}+23m_{18}-9m_{19}+18m_{20}-k,\cr
v_{14} =& 1231m_{15}-672m_{16}+249m_{17}-595m_{18}+259m_{19}-346m_{20}-19k,\cr
v_{15} =& 364m_{16}+19m_{17}+271m_{18}-89m_{19}+290m_{20}-41k,\cr
v_{16} =& 529m_{17}+361m_{18}+185m_{19}+162m_{20}-107k, \cr
v_{17} =& 62m_{18}+m_{19}-22m_{20}+8k, \cr
v_{18} =& 30m_{19}-9m_{20}-8k, \cr
v_{19} =& 3m_{20}-4k.\cr
}
$$
}
%{v[1] = 2*m[2]+m[5]+m[7]+m[10]+m[12]+m[14]+m[15]+m[16]+m[17]+m[18]+2*m[20]-k,
%v[2] = 4*m[3]-m[4]+2*m[5]+2*m[7]+2*m[10]+2*m[12]+2*m[14]+2*m[15]+2*m[16]+m[17]+2*m[18]+2*m[19]+3*m[20]-2*k,
%v[3] = 7*m[4]-2*m[5]+2*m[7]+2*m[10]+2*m[12]+2*m[14]+2*m[15]+2*m[16]+m[17]+2*m[18]+2*m[19]+3*m[20]-2*k,
%v[4] = 33*m[5]-14*m[6]+9*m[7]-14*m[8]+9*m[10]+9*m[12]+9*m[14]+9*m[15]+9*m[16]+15*m[17]+9*m[18]+16*m[19]+24*m[20]-9*k,
%v[5] = 52*m[6]-24*m[7]-14*m[8]+9*m[10]+9*m[12]+9*m[14]+9*m[15]+9*m[16]+15*m[17]+9*m[18]+16*m[19]+24*m[20]-9*k,
%v[6] = 24*m[7]+m[8]+4*m[10]+4*m[12]+4*m[14]+4*m[15]+4*m[16]+11*m[17]-9*m[18]-3*m[19]+2*m[20]-4*k,
%v[7] = 35*m[8]+8*m[10]+8*m[12]+8*m[14]+8*m[15]+8*m[16]+13*m[17]+9*m[18]+15*m[19]+22*m[20]-8*k,
%v[8] = 2*m[9]-m[10],
%v[9] = 211*m[10]-140*m[11]+m[12]+m[14]+m[15]+m[16]+41*m[17]+23*m[18]+50*m[19]+64*m[20]-k,
%v[10] = 282*m[11]-210*m[12]+m[14]+m[15]+m[16]+41*m[17]+23*m[18]+50*m[19]+64*m[20]-k,
%v[11] = 119*m[12]-94*m[13]+m[14]+m[15]+m[16]-53*m[17]+23*m[18]+50*m[19]-30*m[20]-k,
%v[12] = 144*m[13]-118*m[14]+m[15]-118*m[16]-53*m[17]+23*m[18]-69*m[19]-30*m[20]-k,
%v[13] = 86*m[14]-71*m[15]-58*m[16]-5*m[17]+23*m[18]-9*m[19]+18*m[20]-k,
%v[14] = 1231*m[15]-672*m[16]+249*m[17]-595*m[18]+259*m[19]-346*m[20]-19*k,
%v[15] = 364*m[16]+19*m[17]+271*m[18]-89*m[19]+290*m[20]-41*k,
%v[16] = 529*m[17]+361*m[18]+185*m[19]+162*m[20]-107*k,
%v[17] = 62*m[18]+m[19]-22*m[20]+8*k,
%v[18] = 30*m[19]-9*m[20]-8*k,
%v[19] = 3*m[20]-4*k}

\smallskip

\noindent After using (\ref{matmatdegree}) to express $m_1$ in terms of
$m_2, \ldots, m_{20}$, and $k$, we can rewrite the equation $C^2=2g_a-2$ as 
\begin{equation}\label{matDsq}
\eqalign{
112&(3-3g_a+k^2) 
=84v_1^2 +42v_2^2+6v_3^2 +\frac{4v_4^2}{11} + \frac{14v_5^2}{143} +
\frac{7v_6^2}{13} + \cr
&+\frac{v_7^2}{5}+84v_8^2+\frac{6v_9^2}{1055}+\frac{28v_{10}^2}{9917}+
\frac{12v_{11}^2}{799}+\frac{v_{12}^2}{102}+\frac{7v_{13}^2}{258}+\cr
&+\frac{7v_{14}^2}{52933}+\frac{6v_{15}^2}{16003}+\frac{6v_{16}^2}{6877}+
\frac{336v_{17}^2}{16399}+\frac{28v_{18}^2}{155}+\frac{28v_{19}^2}{5}. \cr
}
\end{equation}

\noindent Suppose $k$ and $g_a$ are fixed. Since the $m_i$ are all 
integral, so are the $v_j$. As the right-hand side of
(\ref{matDsq}) is a positive definite quadratic form in the $v_j$, we
find that there are only finitely many integral solutions $(v_1,
\ldots, v_{19})$ of (\ref{matDsq}). The $m_i$ being linear
combinations of the $v_j$, there are also only finitely many integral
solutions in terms of the $m_i$. In our case the even degree $d$ is smaller
than $4$, so $d=2$ and $k=1$. As all curves have even degree, the
conic $C$ is irreducible and hence, as all irreducible conics are, smooth.
Therefore we have $g_a=0$. A computer search shows 
that for $k=1$ and $g_a=0$ there are exactly $441$ solutions
of (\ref{matDsq}) corresponding to integral $m_i$. 

%Any two different effective divisors have a nonnegative intersection
%number, so an effective divisor with negative self-intersection is not
%linearly equivalent to any other effective divisor. We conclude that
%there are at most $441$ effective divisors $D$ on $\bXbar$ satisfying
%$D^2=-2$ and $H \cdot D = 2$.
%On the other hand, we already know $441$ of such divisors. Of 
%those, $9$ correspond to the curves in the orbit of $D_{17}$. 

By Lemma \ref{matcurvesonK3}(a) these correspond to $441$ effective
divisor classes $[D]$ on $\bXbar$ with $D^2=-2$ and $H\cdot D =
2$. We will exhibit $441$ of such divisors satisfying the hypotheses
of Lemma \ref{matcurvesonK3}(b). That lemma then implies that 
each is the only effective 
divisor in its equivalence class and we conclude that they are the
only $441$ effective divisors $D$ on $\bXbar$ satisfying $D^2=-2$ and
$D \cdot H = 2$. 

The first $9$ of these $441$ divisors correspond to 
the curves in the orbit of $D_{17}$. Another $16$ correspond to 
$D_{10}+\varepsilon_1 E_1+\varepsilon_2 E_2+\varepsilon_3 E_3+
\varepsilon_4 E_4$ where $\varepsilon_i 
\in\{0, 1\}$ and the $E_i$ are the four exceptional curves of $\pi$ that 
meet $D_{10}$. Each of these $16$ divisors
generates an orbit under $G$ of size $18$, giving $288$ divisors 
on $\bXbar$ altogether. The last $144$ divisors correspond to
the divisors in the size $36$ orbits of $D_{16} + \delta_1 M_1 + 
\delta_2 M_2$, with 
$\delta_i \in\{0, 1\}$ and where $M_1$ and $M_2$ are the exceptional curves
of $\pi$ in the fiber above $t=2$. Of these $441$ effective divisors,
only $63$ are the strict transform of a curve on $\Xbar$, all in an
orbit of one of the curves described in (\ref{matmatcurves}).
\end{proofof}

\bigskip

\mysectiononumber{\titleref}\label{matsecref}

\begin{itemize}
\setlength{\itemsep}{-1mm}
%\myrefart{\aas}{Aassila, M.}{Some results on Heron
%  triangles}{Elem. Math.}{{\bf 56} (2001)}{143--146}
\myrefart{\artin}{Artin, M.}{On Isolated Rational Singularities of
  Surfaces}{Amer. J. Math.}{{\bf 88} (1966)}{129--136}
\myrefart{\blv}{Beukers, F., van Luijk, R. and Vidunas, R.}{A linear algebra
exercise}{Nieuw Archief voor Wiskunde}{{\bf 3} (2002)}{139--140}
\myrefbook{\BPV}{Barth, W., Peters, C., and Van de Ven, A.}{Compact
  Complex Surfaces}{Ergebnisse der Mathematik und ihrer Grenzgebiete,
  3. Folge, Band 4, Springer-Verlag}{1984}
%\myrefbook{\neron}{Bosch, S., L\"utkebohmert, W., and Raynaud,
%  M.}{N\'eron Models}{Springer-Verlag, Berlin}{1990}
%\myrefart{\bombmum}{Bombieri, E. and Mumford, D.}{Enriques' classification of
%surfaces in char. $p$, II}{Complex Analysis and Algebraic
%Geometry--Collection of papers dedicated to K. Kodaira}{ed. W.L. Baily and T.
%Shioda, Iwanami and Cambridge Univ. Press (1977)}{23--42}
\myrefart{\bremner}{Bremner, A.}{On squares of squares II}{Acta
  Arith.}{{\bf 99}, no. 3 (2001)}{289--308}
%\myrefart{\bruce}{Bruce, J. and Wall, C.}{On the classification of
%  cubic surfaces}{J. London Math. Soc. (2)}{{\bf 19} (1979)}{245--256}
\myrefart{\BT}{Bogomolov, F. and Tschinkel, Yu.}{Density of rational
  points on elliptic K3 surfaces}{Asian J. Math.}{{\bf 4}, 2
  (2000)}{351--368}  
\myrefart{\chin}{Chinburg, T.}{Minimal Models of Curves over Dedekind
rings}{Arithmetic Geometry}{ed. Cornell, G. \& Silverman,
  J. (1986)}{309--326} 
%\myrefart{\deligne}{Deligne, P.}{La Conjecture de
%  Weil. I}{Publ. Math. IHES}{{\bf 43} (1974)}{273--307}
%\myrefart{\duval}{Du Val, P.}{On isolated singularities which do not
%  affect the conditions of adjunction, Part I}{Proc. Cambridge
%  Phil. Soc.}{{\bf 30} (1934)}{453--465}
%\myrefbook{\egaII}{Grothendieck, A.}{\'El\'ements de g\'eom\'etrie
%  alg\'ebrique. IV. \'Etude globale
%  \'el\'ementaire de quelques classes de morphismes}
%{Inst. Hautes \'Etudes Sci. Publ. Math., no. {\bf 8}}{1961}
%\myrefbook{\egaone}{Grothendieck, A.}{\'El\'ements de g\'eom\'etrie
%  alg\'ebrique. IV. 
%\'Etude locale des sch\'emas et des morphismes de sch\'emas, Premi\`ere
%  partie}{Inst. Hautes 
%\'Etudes Sci. Publ. Math., no. {\bf 20}}{1964}
\myrefbook{\egatwo}{Grothendieck, A.}{\'El\'ements de g\'eom\'etrie
  alg\'ebrique. IV.
\'Etude locale des sch\'emas et des morphismes de sch\'emas, Seconde
  partie}{IHES Publ. Math., no. {\bf 24}}{1965}
%\myrefbook{\egafour}{Grothendieck, A.}{\'El\'ements de g\'eom\'etrie
%  alg\'ebrique. IV.
%\'Etude locale des sch\'emas et des morphismes de sch\'emas, Quatri\`eme
%  partie}{Inst. Hautes 
%\'Etudes Sci. Publ. Math., no. {\bf 32}}{1967}
%\myrefbook{\fult}{Fulton, W.}{Intersection  Theory, second
%  edition}{Springer}{1998} 
%\myrefbook{\guy}{Guy, R.K.}{Unsolved Problems in Number
%  Theory}{Problem Books in Math., Springer-Verlag, New-York}{1994}
%\myrefart{\hagtwo}{Hartshorne, R.}{Equivalence relations of algebraic
%  cycles and subvarieties of small codimension}{Algebraic Geometry,
%  Arcata 1974}{Amer. Math. Soc. Proc. Symp. Pure Math. {\bf
%  29} (1975)}{129--164} 
\myrefbook{\hag}{Hartshorne, R.}{Algebraic Geometry}{GTM {\bf 52},
Springer-Verlag, New-York}{1977}
\myrefart{\inose}{Inose, H.}{On certain Kummer surfaces which can be
  realized as non-singular quartic surfaces in $\P^3$}{Journal of the
  Faculty of Science. The 
University of Tokyo}{Section 1A, mathematics, {\bf 23} (1976)}{545--560}
%\myrefart{\kramluc}{Kramer, A.-V. and Luca, F.}{Some remarks on Heron
%  triangles}{Acta Acad. Paedagog. Agriensis Sect. Mat. (N.S.)}{{\bf
%  27}  (2000)}{25--38 (2001)} 
%\myrefart{\kodtwo}{Kodaira, K.}{On compact analytic surfaces II-III}{Ann. of
%Math.}{{\bf 77} (1963), pp. 563--626; {\bf 78} (1963)}{1--40}
%\myrefart{\kod}{Kodaira, K.}{On the structure of compact complex
%  analytic surfaces I, II}{Amer. J. Math.}{{\bf 86} (1964),
%  pp. 751--798; {\bf 88} (1966)}{682--721} 
%\myrefart{\lich}{Lichtenbaum, S.}{Curves over discrete valuation
%  rings}{Amer. J. Math.}{{\bf 90} (1968)}{380--405}
%\myrefart{\lip}{Lipman, J.}{Rational singularities, with applications
%  to algebraic surfaces and unique
%  factorization}{IHES Publ. Math.}{{\bf 36} (1969)}{195--279}
%\myrefbook{\luijk}{van Luijk, R.}{On Perfect Cuboids}{2000}
%{\newline {\tt http://www.math.leidenuniv.nl/reports/2001.html},
%  MI-2001-12}
\myrefbook{\luijkheron}{van Luijk, R.}{An elliptic K3 surface
associated to Heron triangles}{To be published}{2004}
%\myrefbook{\manin}{Manin, Y.}{Cubic forms: Algebra, Geometry,
%  Arithmetic}{North-Holland, Amsterdam}{1974}
%\myrefart{\maz}{Mazur, B.}{Modular curves and the Eisenstein ideal}{IHES Publ.
%Math.}{{\bf 47} (1977)}{33--186}
%\myrefbook{\milne}{Milne, J.S.}{\'Etale Cohomology}{Princeton
%  Mathematical Series 
%  {\bf 33}, Princeton University Press, New Jersey}{1980}
%\myrefart{\nagata}{Nagata, M.}{On rational surfaces I,
%  II}{Mem. coll. Sci. Kyoto (A)}{{\bf 32} (1960), pp. 351--370; {\bf
%    33} (1960)}{271--293}
\myrefart{\naw}{Problem 10}{Problem Section}{Nieuw Archief voor Wiskunde}
{{\bf 1} (2000)}{413--417}
\myrefart{\nik}{Nikulin, V.}{Integral symmetric bilinear forms and
  some of their applications}{Math. USSR Izvestija}{{\bf 14}, 1
  (1980)}{103--167} 
%\myrefart{\pink}{Pinkham, H.}{Singularit\'es Rationelles de
%  Surfaces}{S\'eminaire sur les Singularit\'es des Surfaces}
%{Lect. Notes in Math. {\bf 777},
%ed. M. Demazure, H. Pinkham, and B. Teissier, Springer-Verlag
%  (1980)}{147--172} 
%\myrefbook{\sull}{O'Sullivan, M.}{Classification and Divisor Class Groups of
%  Normal Cubic Surfaces in $\P^3$}{U.C. Berkeley, Ph.D. dissertation
%  (not published)}{1996}
%\myrefbook{\sgaone}{Grothendieck, A.}{Rev\^etements \'etales et Groupe
%Fondamental}{Lect. Notes in Math. {\bf 224}, Springer-Verlag,
%  Heidelberg}{1971} 
%\myrefbook{\sgafourh}{Grothendieck, A. et al.}{Cohomologie \'etale}
%{Lect. Notes in Math. {\bf 569}, Springer-Verlag,
%  Heidelberg}{1977} 
\myrefbook{\sgasix}{Grothendieck, A. et al.}{Th\'eorie des
  Intersections et Th\'eor\`eme de Riemann-Roch}{Lect. Notes in
  Math. {\bf 225}, Springer-Verlag, Heidelberg}{1971} 
%\myrefbook{\shaf}{Shafarevich, I.}{Lectures on Minimal Models and Birational
%Transformations of Two-dimensional Schemes}{Tata Institute, Bombay}{1966}
\myrefart{\shioda}{Shioda, T.}{On the Mordell-Weil Lattices}{Comm. Math. Univ.
Sancti Pauli}{{\bf 39}, 2 (1990)}{211--240}
\myrefart{\SI}{Shioda, T. and Inose, H.}{On singular K3 surfaces}{Complex
Analysis and Algebraic Geometry}{(1977)}{119--136} 
\myrefbook{\silv}{Silverman, J.H.}{The Arithmetic of Elliptic Curves}{GTM {\bf
106}, Springer-Verlag, New-York}{1986}
\myrefbook{\silvtwo}{Silverman, J.H.}{Advanced Topics in the Arithmetic of
Elliptic Curves}{GTM {\bf 151}, Springer-Verlag, New-York}{1994}
\myrefart{\tate}{Tate, J.}{Algorithm for determining the type of a 
singular fiber in an elliptic pencil}
{Modular functions of one variable IV}{Lect. Notes in Math. {\bf 476},
ed. B.J.~Birch and W.~Kuyk, Springer-Verlag, Berlin (1975)}{33--52}
\end{itemize}

\end{document}